\documentclass[11pt]{article}
\usepackage{amsmath,amsfonts,amssymb}

\newtheorem{theorem}{Theorem}[section]
\newtheorem{corollary}[theorem]{Corollary}
\newtheorem{lemma}[theorem]{Lemma}
\newtheorem{sublemma}[theorem]{Sublemma}
\newtheorem{proposition}[theorem]{Proposition}

\newenvironment{proof}[1][Proof]{\vskip -5pt \noindent\textbf{#1.} }{\
\rule{0.5em}{0.5em}}

\begin{document}

\author{Sun-Chin Chu\\Department of Mathematics \\
National Chung Cheng University\\
Chia-Yi, Taiwan\\
scchu@math.ccu.edu.tw}

\date{  }

\title{Type II ancient solutions to the Ricci
flow on surfaces}

\maketitle

\section{Introduction}

Recall that a solution to the Ricci flow is called ancient if it
exists on a time interval $(-\infty, \omega)$ containing $t_0$ for
some $t_0 \in (-\infty,\infty)$.  Let $(\mathcal{M}^n,g(t))$ be a
solution to the Ricci flow. We define what it means to be a Type I
or Type II ancient solution as follows: \medskip

First, $(\mathcal{M}^n,g(t))$ is a complete ancient solution with
bounded curvature ({\it the bound may depend on time}.)
\begin{itemize}
  \item It is Type I if  it satisfies
  \begin{equation}
\sup_{\mathcal{M}\times (-\infty,t_0]} |t| | \mathrm{Rm}(x,t) | <
\infty.
  \end{equation}
  \item It is Type II  if it satisfies
  \begin{equation} \label{eq:ast}
\sup_{\mathcal{M}\times (-\infty,t_0]} |t| | \mathrm{Rm}(x,t) | =
\infty.
  \end{equation}
\end{itemize}
Note that hypothesis \eqref{eq:ast} on a Type II ancient solution
implies that the metric must be non-flat.

In \cite{H_sing}, Hamilton shows that the only Type I ancient
solutions on surfaces are the round sphere $\mathbb{S}^2$ and the
flat plane $\mathbb{R}^2$, and their quotients. Therefore, any
non-flat complete ancient solution with bounded curvature on a
noncompact surface is Type II. That is, there does not exist a Type
I non-flat ancient solution on a noncompact surface. As we know so
far, Type II ancient solutions on surfaces have not been classified.
It is conjectured that the noncompact case should correspond to the
cigar soliton and compact case to the Rosenau solution \cite{R}.

Of particular interest to us is to study Type II ancient solutions
on surfaces. By the strong maximum principle, we see that
$R(g(t))\equiv 0$ everywhere and it is flat, or $R(g(t))>0$
everywhere and it is diffeomorphic to $\mathbb{S}^2$ or
$\mathbb{R}^2$. Here and throughout, let $(\mathbb{R}^2,g(t))$
denote a Type II ancient solution to the Ricci flow. We shall see
that such a solution can be extended to a complete eternal solution,
i.e., it is defined on $(-\infty,\infty)$. Remark that the curvature
is still bounded at each time slice.

The paper is organized as follows. In section 2, we study the limits
backwards in time, in a way analogous to a maximal solution of Type
IIb in \cite{H_sing}, of Type II ancient solutions on surfaces.
Proposition \ref{thm:backward_limit=cigar} shows that the backward
limit of such a solution is a multiple of the cigar soliton. In
section 3, we investigate the asymptotic volume ratio, total
curvature, aperture and circumference at spatial infinity of Type II
ancient solutions on complete noncompact surfaces. We shall see that
the scalar curvature of such a solution decays to zero at spatial
infinity, hence that these quantities are preserved under the Ricci
flow. Theorem \ref{thm:zero_aperture=finite_circumference} shows
that the circumference at spatial infinity of such a solution is
finite, therefore, the volume grows linearly. Since Riemann surfaces
are K\"{a}hler, this improves Ni's theorem \cite{Ni}, namely that
any non-flat ancient solution to the K\"{a}hler-Ricci flow with
bounded nonnegative bisectional curvature has asymptotic volume
ratio zero.  By the Harnack estimate,  the function
$R_{\max}(t)=\sup R(\cdot,t)$ is nondecreasing. Does a Type II
ancient solution on a surface satisfy $\lim_{t\rightarrow
-\infty}R_{\max}(t)>0$? Theorem \ref{thm:c=2pi=>r>=4} gives an
affirmative answer to the noncompact case.

\section{Taking limits backwards in
time}\label{sect:backward_limit}

In this section, we shall take limits backwards in time of Type II
ancient solutions on surfaces.

\subsection{The compactness theorem} \label{subsect:compactness_theorem}

Recall the definition of convergence of pointed solutions to the
Ricci flow.

To begin with, we fix a time interval $(\alpha,\omega)$ with
$-\infty \leq \alpha <0$ and $0<\omega \leq \infty$.

\medskip

\noindent\textbf{Definition}. A sequence
$\{(\mathcal{M}^n_k,g_k(t),O_k)\}_{k\in \mathbb{N}},\ t\in
(\alpha,\omega)$, of complete pointed solutions to the Ricci flow
converges to a complete pointed solution to the Ricci flow
$(\mathcal{M}^n_\infty,g_\infty(t),O_\infty),\ t\in (\alpha,\omega)$,
if there exist

\begin{itemize}

\item[(1)] an exhaustion $\{U_k\}_{k\in \mathbb{N}}$ of $\mathcal{M}_\infty$
by open sets with $O_\infty \in U_k$, and

\item[(2)] a sequence of diffeomorphisms $\Phi_k:U_k \rightarrow
V_k=\Phi_k(U_k) \subset \mathcal{M}_k$ with $\Phi_k(O_\infty)=O_k$
such that $(U_k,\Phi_k^\ast[g_k(t)|_{V_k}])$ converges in $C^\infty$
to $(\mathcal{M}_\infty,g_\infty(t))$ uniformly on compact sets in
$\mathcal{M}_\infty \times (\alpha,\omega)$.
\end{itemize}

We review Hamilton's compactness theorem for sequences of solutions
to the Ricci flow as follows.

\bigskip

\noindent\textbf{Theorem} (Hamilton \cite{H_sing}). \textit{Let}
$\{(\mathcal{M}^n_k,g_k(t),O_k)\}_{k\in \mathbb{N}},\ t \in
(\alpha,\omega)$, \textit{be a sequence of complete pointed
solutions to the Ricci flow such that}

\begin{itemize}
\item[(i)] (\textit{uniformly bounded curvatures})
\[
|\mathrm{Rm}_k|_k \leq C_0 \textit{ on } \mathcal{M}_k \times
(\alpha,\omega)
\]
\textit{for some constant} $C_0<\infty$ \textit{independent of}
$k$, \textit{and}

\item[(ii)] (\textit{injectivity radius estimate at} $t=0$)
\[
\mathrm{inj}(O_k,g_k(0)) \geq i_0
\]
\textit{for some constant} $i_0 >0$.
\end{itemize}

\textit{Then there exists a subsequence} $\{j_k\}_{k \in
\mathbb{N}}$ \textit{such that} $\{
(\mathcal{M}_{j_k},g_{j_k},O_{j_k})\}$ \textit{converges to a
complete pointed solution to the Ricci flow}
$(\mathcal{M}^n_\infty,g_\infty(t),O_\infty),$ $ t\in
(\alpha,\omega)$, \textit{as} $k \rightarrow \infty$.

\medskip

\noindent\textbf{Remark}. (1) In fact, if there is a subsequence
$(\mathcal{M}_{j_k},g_{j_k}(0),O_{j_k})$ convergent to a limit
$(\mathcal{M}_\infty,g_\infty,O_\infty)$, then there is a
subsequence which converges at all times. (2) For the Ricci flow it is known that curvature bounds on $(\alpha,\omega)$ imply bounds on all derivatives of the curvature on $[\alpha + \varepsilon,\omega)$ for any $\varepsilon>0$. Thus we need only assume the
curvature bound for solutions to the Ricci flow.

\subsection{The backward limit} \label{subsect:rosenau}

Let's first take a look at the Rosenau solution \cite{R}.

Let $\mathcal{M}^2$
be the cylinder $\mathbb{R}\times \mathbb{S}_1$, where
$\mathbb{S}_1$ is the circle of radius 1. We define a  solution
$g(x,\theta,t)$, $t<0$, to the Ricci flow on $\mathcal{M}^2$  by
\[
g(x,\theta,t)=\frac{\sinh(-t)}{\cosh x + \cosh t}(dx^2+d\theta^2).
\]
It is easy to justify that the solution  $g(x,\theta,t)$ extends to
a complete ancient solution to the Ricci flow on the sphere
$\mathbb{S}^2$. This complete ancient solution on $\mathbb{S}^2$ is
the so-called Rosenau solution.

By
straightforward computation, the scalar curvature of the metric is given
by
\[
R(x,\theta,t)=\frac{1+\cosh t \cdot \cosh x }{\sinh (-t)(\cosh x+
\cosh t)}>0,
\]
and attains its maximum curvature at the poles $x=\pm \infty$:
\begin{align*}
R_{\max}(t) & =\lim_{|x| \rightarrow \infty} \frac{1+\cosh t
\cdot \cosh x }{\sinh (-t)(\cosh x+ \cosh t)} \\
& =\coth (-t) >0
\end{align*}
for all $t<0$.  Since $\lim_{t\rightarrow 0^{-}}R_{\max}(t)=\infty$,
the Rosenau solution is ancient, but not eternal. Note that the Rosenau solution has a Type I singularity as $t\nearrow 0$. By the
fact that
\[
\lim_{t\rightarrow 0^{-}}\frac{R(x,\theta,t)}{R_{\max}(t)}=1
\text{ for all } (x,\theta)\in \mathbb{R}\times \mathbb{S}_1,
\]
the normalized  solution converges to the round sphere
$\mathbb{S}^2$ as $t \rightarrow 0^{-}$.  On the other hand, we have
\[
\sup_{\mathbb{S}^2 \times (-\infty,-1]}|t|R
=\sup_{(-\infty,-1]}|t|\coth(-t)=\infty,
\]
which means that it is a   Type II ancient solution on $\mathbb{S}^2$.

\medskip

To study the limits  backwards in time, in a way analogous to a
maximal solution of Type IIb in \cite{H_sing}, of Type II ancient
solutions  on surfaces, we need the following.

\begin{lemma} \label{lemma:backward_limit}
Suppose that $(\mathcal{M}^n,g(t))$ is a Type II ancient solution
and satisfies the injectivity radius bound, namely
\begin{equation} \label{eq:star}
\mathrm{inj}(\mathcal{M},g(t)) \geq \frac{c}{\sqrt{\mathcal{K}(t)}}
\ \ \ \ \text{ at each time } t,
\end{equation}
where $c$ is a positive constant and
$$
\mathcal{K}(t)=\sup_{x\in \mathcal{M}}  |
\mathrm{Rm}(x,t) |.
$$
Then there exists a sequence of dilations of the solution which
converges to a Type II singularity model.
\end{lemma}

\begin{proof}
Let $\{ \gamma_j \}$ be any sequence with $\gamma_j \nearrow 1$, and
choose any sequence of time $T_{j}\searrow-\infty$, and pick
$(x_{j}, t_{j}) \in \mathcal{M}\times [T_i,0]$ such that
\[
|t_{j}|(t_{j}-T_{j})|\mathrm{Rm}(x_{j},t_{j})|\geq \gamma_{j}
\sup\limits_{M\times [T_{j},0]}\ |t|(t-T_{j})|\mathrm{Rm}(x,t)|.
\]

Now consider the dilated solutions
\begin{equation*}
 g_{j}(t)=| \mathrm{Rm}(
x_{j},t_{j}) | \cdot g\left( t_{j}+\frac{t}{| \mathrm{Rm}(
x_{j},t_{j}) |} \right).
\end{equation*}
By the injectivity radius bound \eqref{eq:star} and definition of
$g_j$, we obtain
\[
\mathrm{inj}(x_j,g_j(0)) \geq c.
\]
Each solution $g_{j}$ exists on the  time interval $(-\infty
,\frac{ \omega -t_{j}}{ | \mathrm{Rm}( x_{j},t_{j}) |} )$, which
contains the subinterval $[-\alpha _{j},\omega _{j}]$ with
\[
\alpha_{j}  = \left( t_{j}-T_{j}\right) \left| \mathrm{Rm}\left(
x_{j},t_{j}\right) \right| \text{ and } \omega_{j}  = -t_{j}\left|
\mathrm{Rm}\left( x_{j},t_{j}\right) \right|.
\]
By definition, we have
\begin{align*}
\frac{1}{1/\alpha_{j}+1/\omega_{j}} & =\frac{\alpha_{j}\omega_{j}}{%
\alpha_{j}+\omega_{j}}=\frac{\left| t_{j}\right| \left(
t_{j}-T_{j}\right) \left| \mathrm{Rm}\left( x_{j},t_{j}\right)
\right| }{\left|
T_{j}\right| } \\
& \geq\gamma_{j} \sup_{M\times%
\lbrack T_{j},0]}\frac{\left| t\right| \left( t-T_{j}\right) \left|
\mathrm{Rm}\left( x,t\right) \right| }{\left| T_{j}\right| } \\
& \rightarrow\infty \ \ \ \ \text{ as }j\rightarrow\infty \text{
since } T_{j}\rightarrow-\infty.
\end{align*}
This implies that $\lim_{j\rightarrow \infty} \alpha_j=\infty$ and
$\lim_{j\rightarrow \infty} \omega_j=\infty$. Therefore, for any
given $\alpha,\, \omega>0$, the interval $(-\alpha_j,\omega_j)$
contains the subinterval $(-\alpha,\omega)$ for $j$ sufficiently
large.

On the other hand, for all $(x,t) \in \mathcal{M}\times (
-\alpha,\omega)$, we see that
\begin{align*}
| \mathrm{Rm}_{j}( x,t) | & =\frac{1}{| \mathrm{Rm}( x_{j},t_{j})|
}\left| \mathrm{Rm}\left( x,t_{j}+\frac{t}{| \mathrm{Rm}
( x_{j},t_{j}) | }\right) \right|  \\
& \leq \frac{( t_{j}-T_{j}) | t_{j}| }{\gamma_j \left(
t_{j}-T_{j}+\frac{t}{| \mathrm{Rm}( x_{j},t_{j}) |} \right) \left| t_{j}+
\frac{t}{| \mathrm{Rm}( x_{j},t_{j}) |}\right | } \\
& =\frac{\alpha _{j}\omega _{j}}{\gamma_{j} ( \alpha _{j}+t) (
\omega _{j}-t) }
\end{align*}
uniformly bounded for $j$ sufficiently large since $\alpha
_{j}\rightarrow \infty $, $\omega _{j}\rightarrow \infty $, and
$\gamma _{j}\nearrow 1$.

Consequently, we conclude that the sequence
$(\mathcal{M},g_j(t),x_j)$, $-\alpha<t<\omega$, satisfies the
hypotheses of Hamilton's compactness theorem. It follows that there
exits a subsequence of $(\mathcal{M},g_j(t),x_j)$ which limits to a
complete pointed eternal solution $(\bar{\mathcal{M}},\bar{g}(t)
,\bar{x}) $ satisfying
\begin{equation*}
\sup_{\bar{\mathcal{M}}\times ( -\infty ,\infty ) }| \mathrm{Rm}|
\leq 1=| \mathrm{Rm}( \bar{x},0)| ,
\end{equation*}
that is, the limit is a Type II singularity model. The lemma
follows. \hfill
\end{proof}

\medskip

Since a Type II  ancient solution on a surface has positive
curvature everywhere, the metric satisfies the injectivity radius
bound \eqref{eq:star}.  Lemma \ref{lemma:backward_limit} implies
that the backward limit of such a solution is a Type II singularity
model. By construction, the only curvature is positive and attains
its maximum in space-time, therefore, it follows from \cite{H_sing}
that the limit must be a multiple of the cigar soliton. We conclude
this section with the following.

\begin{proposition} \label{thm:backward_limit=cigar}
If a complete  ancient solution to the Ricci flow on a surface  with
bounded curvature is not a quotient of the round sphere or of the flat
plane, then the ancient solution is   Type II. Moreover, the backward
limit of such a solution is a multiple of the cigar soliton.
\end{proposition}

As a corollary,  we see that  the backward limit of the Rosenau
solution is the cigar soliton.

\medskip

\noindent \textbf{Remark}. Proposition
\ref{thm:backward_limit=cigar} is also obtained independently by Chow, Lu
and Ni
\cite{CLN}.

\section{The geometry at spatial infinity of Type II ancient solutions
on $\mathbb{R}^2$}\label{sect:geometry_at_infty}

\setcounter{equation}{0} In this section, first we recall the
asymptotic volume ratio,   total curvature, aperture,  and
circumference at infinity of complete noncompact surfaces with
bounded positive curvature. Next, we  study these quantities      of
Type II  ancient solutions on $\mathbb{R}^2$.   We shall see that
these quantities are preserved under the Ricci flow. L. Ni \cite{Ni}
proves that any non-flat ancient solution to the K\"{a}hler-Ricci
flow with bounded nonnegative bisectional curvature has asymptotic
volume ratio zero. For Riemann surfaces, Theorem
\ref{thm:zero_aperture=finite_circumference} improves Ni's theorem
since finite circumference at infinity implies that  the volume
grows linearly, hence that the asymptotic volume ratio is zero.

\subsection{The geometry of complete  surfaces at infinity}
\label{subsect:aperture_circumeference_total_curvature}

Suppose that $(\mathcal{M}^n,g)$ is a complete Riemannian manifold
with nonnegative Ricci tensor. The Bishop-Gromov theorem  says that
the function
\[
r\rightarrow\frac{\mathrm{Vol}(B(p,r))}{\omega_{n}r^{n}},
\]
where $B(p,r)=\{ x \ |\ \mathrm{dist}(x,p) < r \}$, is monotone
decreasing for any $p\in\mathcal{M}$. The asymptotic volume ratio
$\alpha_{g}$ is defined by
\[
\alpha_{g}=\lim_{r\rightarrow\infty}\frac{\mathrm{Vol}
(B(p,r))}{\omega_{n}r^{n}},
\]
which is independent of $p$ and invariant under dilation. It is
known that
\begin{equation}  \label{vol-lower-upper}
\alpha_{g}\omega_{n}r^{n}\leq\mathrm{Vol}(B(p,r))\leq\omega_{n}
r^{n}.
\end{equation}

Suppose that $(\mathbb{R}^2,g)$ is a complete surface with bounded
positive curvature. Let $o\in \mathbb{R}^2$ be some  point
which we call the origin. Denote by $B_s$  the open ball of radius
$s$ around the origin $o$, $\ell(s)$ the length of $\partial B_s$
and $A(s)$ the area of $B_s$. Recall that the total curvature
$\tau_g$ and aperture $\mathcal{A}_g$ of the metric $g$ are given by
\[
\tau_g=\int_{\mathbb{R}^2}Kd\mu_{g} \text{ and
}\mathcal{A}_g=\lim_{s\rightarrow\infty}\frac{\ell(s)}{s},
\]
respectively.   Note that the aperture is also independent of the choice of
the origin and invariant under dilation. It follows from the
Hartman theorem \cite{Hartman} that we have
\begin{subequations}
\begin{equation} \label{eq:t_1}
\lim_{s\rightarrow\infty} \frac{\ell(s)}{s}   =
2\pi-\int_{\mathbb{R}^2}Kd\mu_{g} 
\end{equation}
and
\begin{equation} \label{eq:t_2}
\lim_{s\rightarrow\infty}\frac{A(s)}{2s^2}  =
2\pi-\int_{\mathbb{R}^2}Kd\mu_{g}. 
\end{equation}
\end{subequations}
By the Cohn-Vossen theorem,    the right hand side of
\eqref{eq:t_1} is nonnegative for a complete noncompact convex
surface, that is, the total curvature  is at most $2\pi$. Clearly,
the left hand side of \eqref{eq:t_2} is a multiple of the asymptotic
volume ratio.

From \eqref{eq:t_1} and \eqref{eq:t_2}, for complete noncompact
convex surfaces we see that the aperture is positive if and only if
the asymptotic volume ratio is positive. Since Riemann surfaces are
K\"{a}hler, for a complete non-flat ancient solution with bounded
curvature on $\mathbb{R}^2$, combining \eqref{eq:t_1} and
\eqref{eq:t_2} to Ni's theorem  \cite{Ni} shows that the aperture is
also zero. Therefore, we have that if  the aperture of a complete
ancient solution with bounded curvature on  $\mathbb{R}^2$ is
positive, then the metric is  flat.

Since the scalar curvature of a Type II ancient solution on
$\mathbb{R}^2$ is positive and bounded, by the Bernstein-Bando-Shi
estimates, injectivity radius estimate and the fact that the total
curvature is at most $2\pi$, the scalar curvature of such solutions
decays to zero at spatial infinity. It follows from \cite{H_sing}
that   the  aperture and  asymptotic volume ratio are preserved
under the Ricci flow.  Consequently, the total scalar curvature is
also preserved under the flow.

Now recall that the circumference at infinity of a complete
noncompact surface $(\mathcal{M}^2,g)$  is defined by
\[
\mathcal{C}_{g}=\sup_{K}\inf_{U}\{\ell(\partial U)|\text{
}\forall\text{ compact set }K\subset \mathcal{M},\forall\text{ open
set }U\supset K\}.
\]

For any monotone sequence of compact sets $K_{n}$ exhausting a
complete noncompact surface, we see that
$$
\mathcal{C}_{g}\geq\sup_{n}\inf_{U\supset K_{n}}\ell (
\partial U ).
$$
On the other hand, for any compact set $K$
there exists a $K_{n}$ with $K\subset K_{n}$ so that  we have
$$
\inf_{U\supset K}\ell (\partial U ) \leq\inf_{U\supset K_{n}}\ell(
\partial U ).
$$
Hence it follows that
\begin{equation} \label{eq:C_g=}
\mathcal{C}_{g}=\sup_{n}\inf_{U\supset K_{n}}\ell (
\partial U )  =\lim_{n\rightarrow\infty}\inf_{U\supset
K_{n}}\ell (\partial U).
\end{equation}

Since  the scalar curvature of a Type II ancient solution
on $\mathbb{R}^2$ vanishes at spatial infinity, for any time
interval $[a,b]$ containing $t_0$ there exists a monotone exhaustion
sequence of compact sets $K_{n}$ with
$$
0<R(x,t)<\frac{1}{n} \text{ for } (x,t)\in( \mathbb{R}^{2}\backslash
K_{n})\times[a,b].
$$
Let $\gamma(s)$ be a fixed parameterized curve on $\mathbb{R}^{2}$.
Then the length evolves by the formula%
\[
\frac{d}{dt}\ell_{g(t)} (\gamma )  =\frac{d}{dt}\int_{\gamma}%
\sqrt{g(\partial_s,\partial_s)}ds=\int_{\gamma}-\frac{R}%
{2}\sqrt{g(\partial_s,\partial_s)}ds\text{.}%
\]
This implies that
\[
-\frac{1}{2n}\ell_{g(t)} ( \partial U ) \leq\frac{d}{dt}\ell_{g(t)}(
\partial U )  <0
\]
for any set $U$ with $K_{n}\subset U$ and $t\in [a,b].$ Therefore,
we
have%
\[
e^{-\frac{1}{2n}(t-a)}\inf_{U\supset K_{n}}\ell_{g(a)}( \partial U)
\leq\inf_{U\supset K_{n}}\ell_{t}(  \partial U)  \leq e^{\frac{1}%
{2n}(b-t)}\inf_{U\supset K_{n}}\ell_{g(b)}( \partial U )
\]
for $t\in [a,b]$.

We conclude this with the following.
\begin{lemma} \label{lemma:circumference=2pi}
If the circumference at infinity of an ancient solution (with
bounded curvature at each time slice) on $\mathbb{R}^2$ is finite for some
$t_0$, then it is constant in time.
\end{lemma}

To explore the aperture and circumference at infinity of  Type II ancient
solutions on  $\mathbb{R}^2$, we employ the theory of isometric
embedding to the surface $(\mathbb{R}^2,g(t))$ as follows.

In \cite{P}, Pogorelov  shows that every complete smooth metric with
positive curvature, given on a plane, is realizable as an unbounded
smooth convex surface $\mathcal{M}^2$ in $\mathbb{R}^3$. By a result
of Stoker \cite{St},   coordinates in $\mathbb{R}^{3}$ can be so
chosen that $\{ z =0 \}\equiv\mathbb{R}^2$ is a supporting
hyperplane to $\mathcal{M}$ at the origin $O=(0,0,0)\in
\mathbb{R}^3$, and $\mathcal{M}$ is the graph of a nonnegative
strictly convex function $f:\Omega \subset \{z=0 \} \rightarrow
\mathbb{R},$ where $\Omega$ is the image of $\mathcal{M}$ under the
orthogonal projection $\pi:\mathbb{R}^{3}\rightarrow \{ z=0 \}.$ Let
$o\in \mathbb{R}^2$ denote the preimage of the origin $O$ under the
map $\mathcal{I}$. Thus, we can identify the point $o$ with the
origin $O$. In what follows, we shall freely realize without
explicit mention a pointed surface $(\mathbb{R}^2,g,o)$ as the graph
of a nonnegative strictly convex smooth function $f$ as above.

Now take the graph of $f$ over the sublevel set $\{f \leq n \}$ as
the compact set $K_n$. By a result of Greene and Shiohama \cite{GS},
the length of level sets is monotone increasing. Together with
\eqref{eq:C_g=} and the observation that $ \inf_{U\supset K_{n}}\ell
( \partial U ) = \ell ( \{ f=n \} ), $ this implies that
$\mathcal{C}_g=\ell(\partial\Omega)$.

There is an essential difference between surfaces with $\tau=2\pi$
and surfaces with $\tau<2\pi$. Considering a convex cone as an
example, it is known that a complete metric with positive curvature
given on a plane may be realized by unbounded convex surface in more
than one way. This is always the case  \cite{O} if the total
curvature of the manifold is less than $2\pi$. On the other hand, it
is known \cite{P}  that for any complete noncompact surface of
nonnegative curvature with $\tau=2\pi$, there is a unique complete
convex surface in $\mathbb{R}^3$ isometric to it up to congruence.
In particular, if $\mathcal{C}_g <\infty$, say $\mathcal{C}_g=2\pi,$
then the embedded surface is inside a circular cylinder of radius
$\pi$. This implies that the tangent cone of the surface is a ray,
thus we see that $\mathcal{A}_g=0$ and $ \tau_{g}=2\pi$ by
\eqref{eq:t_1}. Therefore, we have  that  the embedding is always
rigid if we have $\mathcal{C}_{g(t)}<\infty$.

\medskip

As a corollary of Lemma \ref{lemma:circumference=2pi} , we have the following.
\begin{lemma} \label{prop:embbedding_rigid}
The isometric embedding of a Type II ancient solution $(\mathbb{R}^2,g(t))$ is
 rigid for all $t$  provided that we have $\mathcal{C}_{g(t_0)}< \infty$
 for some $t_0\in (-\infty,\omega)$.
\end{lemma}

\medskip

\noindent\textbf{Remark}. In \cite{DH}, Daskalopoulos and Hamilton
introduce the width $w(g)$ of a metric $g$ on the plane. Let
$F:\mathbb{R}^2 \rightarrow [0,\infty)$ denote a proper function,
i.e., $F^{-1}(c)$ is compact for every $c\in  [0,\infty)$. The width
of $F$ is given by
$$
w(F) = \sup_c \ \ell(\{F = c\}).
$$
Then, the width $w(g)$ is given by the infimum of $w(F)$ over all
smooth proper functions $F$, i.e.,
$$ w = \inf_F w(F).$$ It is clear that if the metric is complete and
has positive curvature, then we have
$$
w(g)=\mathcal{C}_g
$$
since the surface is realizable as the graph of a (proper) strictly convex
function.

\subsection{Circumference of $(\mathbb{R}^2,g(t))$ at infinity}
\label{subsect:zero_aperture}

\medskip

For Riemann surfaces, Theorem
\ref{thm:zero_aperture=finite_circumference}  improves Ni's theorem
 since  the
finiteness of circumference at infinity implies that the volume
grows linearly. Consequently, we  see that  the asymptotic volume ratio
is $0$.

\begin{theorem} \label{thm:zero_aperture=finite_circumference}
The circumference at spatial infinity of a Type II ancient  solution
$(\mathbb{R}^2,g(t))$ is finite and independent of time.
\end{theorem}

\begin{proof} Let's begin the proof with the following observation.

\begin{lemma} \label{lemma:average_upper_bound}
Suppose that $(\mathbb{R}^2,g_0)$ is a complete surface with
positive and bounded curvature. Then there exists a positive
constant $C$ independent of $x\in \mathbb{R}^2$ and $r>0$ such that
\begin{equation}\label{eq:ave}
\frac{1}{\mathrm{Vol}(B(x,r))}\int_{B(x,r)}R\ d\mu \leq \frac{C}{r}.
\end{equation}

\end{lemma}

\begin{proof}
We may assume without loss of generality that $R(x) \leq 2$
for all $x$.  By the injectivity radius estimate of Meyer and Gromoll,
$\text{inj}(\mathbb{R}^2,g_0)$ has a lower bound $\pi$. Therefore,
it follows from Yau's theorem that there exists a  positive constant
$C$ (independent of $x$) such that
\[
\mathrm{Vol}(B(x,r)) \geq Cr \ \ \ \ \text{for } r \geq 1.
\]
Combining this estimate with the Cohn-Vossen theorem shows that
\begin{align*}
 \frac{1}{\mathrm{Vol}(B(x,r))}\int_{B(x,r)}Rd\mu &
 \leq \frac{1}{\mathrm{Vol}(B(x,r))}\int_{\mathbb{R}^2}Rd\mu
\\
& \leq \frac{4\pi}{Cr} \ \ \ \  \text{ for } r \geq 1.
\end{align*}
On the other hand, for $r<1$, it is easy to see that
\[
 \frac{1}{\text{Vol}(B(x,r))}\int_{B(x,r)}Rd\mu
 \leq \frac{1}{\text{Vol}(B(x,r))}\int_{B(x,r)}2d\mu  =2.
\]
Therefore, the lemma follows. \hfill
\end{proof}

\medskip

\noindent\textbf{Remark}. Since Riemann surfaces are K\"{a}hler, it follows from Shi's theorem
\cite{Sh} that the ancient solution $g(t)$ can be extended to an eternal
solution still with bounded curvature at each  time slice. This fact plays a roll in the proof of Lemma \ref{lemma:lower_bound:NT3}.

\medskip

For convenience, let $k(x,r)$ denote the average of the scalar
curvature over $B(x,r)$ with respect to $g(0)$, that is,
\[
k(x,r)=\frac{1}{\mathrm{Vol}(B(x,r))}\int_{B(x,r)}R_{g(0)}
d\mu_{g(0)}.
\]

From the fact that the aperture  of a Type II
ancient solution   on $\mathbb{R}^2$  is zero,
we have $\tau=2\pi$ by \eqref{eq:t_2}, therefore, the isometric
embedding of the surface $(\mathbb{R}^2,g(t))$ is rigid. Let $o$ and
$f$ as given in subsection
\ref{subsect:aperture_circumeference_total_curvature}. Denote the
sublevel set $\{ f \leq n \}$ by $\Omega$, and let $\rho=\text{
dist}(o,\partial \Omega)$ and  $r=\max_{p\in
\partial \Omega}\mathrm{dist}(o,p)$. It is clear that
$\Omega \hookrightarrow (\mathbb{R}^2,g,o)$ is a
compact domain with non-empty convex boundary $\partial \Omega$. We
may combine the results of Lemma 35.3.1 and Theorem 35.3.2 in
\cite{BZ} to conclude that $\Omega$ has a boundary starlike with
respect to $o$ and satisfies
\[
\mathrm{Vol}(\Omega) \geq \frac{\rho}{2} \ell(\partial \Omega).
\]
Combining this estimate with the facts that
$$
\rho \leq r \leq \rho +\ell(\partial \Omega)  \text{ and }
\ell(\partial \Omega)=o(r) \text{ as } r \rightarrow \infty,
$$
we
obtain
\begin{align*}
\mathrm{Vol}(B(o,r)) & \geq \mathrm{Vol}(\Omega)  \geq \frac{\rho}{2}
\ell(\partial \Omega) \\
& > \frac{r}{4} \ell(\partial \Omega)
\end{align*}
for $r$ sufficiently large, and hence
\begin{align*}
 r\cdot \frac{1}{\mathrm{Vol}(B(o,r))}\int_{
B(o,r)} R\ d\mu & < r\cdot \frac{4}{ r \ell(\partial \Omega)}
\int_{ \mathbb{R}^2} R\ d\mu  \\
& = \frac{8\pi}{\ell(\partial \Omega)}.
\end{align*}

Therefore, we obtain the following.
\begin{lemma} \label{lemma:upper_bound}
Let $\Gamma$ be the level curve $\{ f=n \}$, and $r=\max_{p\in
\Gamma}\mathrm{dist}(p,o)$. Then for $r$ sufficiently large, we have
\begin{equation} \label{eq:upper_bound}
rk(o,r) < \frac{8\pi}{\ell( \Gamma)}.
\end{equation}
\end{lemma}

Based on an observation of Ni and Tam (Proposition 2 in \cite{NT3}),
we have the following.

\begin{lemma} \label{lemma:lower_bound:NT3}
\begin{equation} \label{eq:positive_lower_bound}
\limsup_{r\rightarrow \infty} \ rk(o,r)>0.
\end{equation}
\end{lemma}
\medskip

\begin{proof}  The proof proceeds along the same lines as  in  \cite{NT3}.
For the convenience of the reader, we give the proof here.

To obtain a contradiction, we suppose that
\[
\limsup_{r\rightarrow\infty} rk(o,r)=0,
\]
that is,
\[
k(o,r)=o(1/r) \ \ \ \ \text{ as } r\rightarrow \infty.
\]
We shall see that this, in conjunction with Lemma
\ref{lemma:average_upper_bound}, suffices  to claim that the surface
must be flat, in contradiction with hypothesis \eqref{eq:ast} on a
Type II ancient  solution.

Let $(\mathcal{M}^m,g_{\alpha\bar{\beta}}(x,t))$ be a solution to
the K\"{a}hler-Ricci flow. Denote by $F(x,t)$ the log of the volume
element, that is,
\[
F(x,t)=\ \log  \left(\frac{\det (g_{\alpha\bar{\beta}}(x,t))}{\det
(g_{\alpha\bar{\beta}}(x,0))}\right).
\]
We thus have
\begin{equation} \label{eq:F=-Int_R}
F(x,t)=-\int_0^t R(x,\tau)d\tau.
\end{equation}
For convenience, let $\mathfrak{m}(t)=\inf_{\mathcal{M}}
F(\cdot,t)$.

\begin{sublemma}[Ni-Tam \cite{NT3}]  \label{lemma:Prop_1_NT3}
Suppose $(\mathcal{M}^n,g_{\alpha\bar{\beta}})$ is a complete
noncompact K\"{a}hler manifold with nonnegative and bounded
bisectional curvature, and the average function $k$ satisfies
estimate
\begin{equation}\label{eq:ave2}
k(x,r) \leq \frac{C}{r}.
\end{equation}
If there exists some point $x_0\in
\mathcal{M}$ such that
\begin{equation} \label{eq:l_o}
k(x_0,r)=o(1/r) \text{ as } r\rightarrow \infty,
\end{equation}
then we have
\begin{equation} \label{eq:3}
\lim_{t\rightarrow \infty}\frac{-F(x,t)}{t}=0 
\end{equation}
and
\begin{equation} \label{eq:4}
\lim_{t\rightarrow \infty}R(x,t)=0 
\end{equation}
for all $x\in \mathcal{M}$.
\end{sublemma}

\medskip

\begin{proof}
Recall that Shi's theorem
\cite{Sh} implies that the ancient solution $g(t)$ can be extended to an eternal
solution. Since  $k(x_0,r)=o(1/r)$ implies $k(x,r)=o(1/r)$ for all $x\in
\mathcal{M}$, it suffices to prove that estimates $\eqref{eq:3}$ and
$\eqref{eq:4}$ are valid for $x_0$.

It follows from Theorem 7.10 in \cite{Sh}, Corollary 2.1 in
\cite{NT1} and estimate \eqref{eq:ave2} that we have
\[
-\mathfrak{m}(t) \leq C
t^{\frac{1}{2}}(1-\mathfrak{m}(t))^{\frac{1}{2}},
\]
hence that
\begin{equation}  \label{eq:m(t)_upper_bound}
1-\mathfrak{m}(t) \leq C(1+t) \ \ \ \text{ for all } t. 
\end{equation}
This, together with Theorem 2.1 in \cite{NT1}, implies that
\begin{align*}
-F(x_0,t) & \leq
C\left[\left(1+\frac{t(1-\mathfrak{m}(t))}{r^2}\right)\int_0^rsk(x_0,s)ds
-\frac{t\mathfrak{m}(t)(1-\mathfrak{m}(t))}{r^2}\right]
\\
& \leq C \left[\left(1+\frac{t^2}{r^2}\right)\int_0^r sk(x_0,s)ds
+\frac{t^3}{r^2}\right]
\end{align*}
for some constant $C$. By $\eqref{eq:l_o}$, for any given
$\varepsilon>0$, there exists a positive constant $r_0$ such that
$k(x_0,r) \leq \varepsilon /r$ whenever $r>r_0$. Putting $r=t/
\sqrt{\varepsilon}$ in the above inequality, we get
\begin{align*}
-F(x_0,t) & \leq C\left[(1+\varepsilon)\int_0^{r_0}
sk(x_0,s)ds+(1+\varepsilon)\varepsilon(r-r_0)+\varepsilon t\right]
\\
& \leq C \left[\int_0^{r_0} sk(x_0,s)ds
+\sqrt{\varepsilon}t+\varepsilon t\right]
\end{align*}
for $t$ sufficiently large. By dividing both sides of the inequality
by $t$, and letting  $t\rightarrow \infty$ and then $\varepsilon
\rightarrow 0$, we conclude that estimate $\eqref{eq:3}$ holds for
$x_0$.

The trace Harnack inequality says that the function $t\mapsto
tR(x,t)$ is increasing in time, hence that
\[
sR(x,s)\cdot\frac{1}{t} \leq R(x,t) \ \text{ for } 0<s<t.
\]
Integrate over $t$ from 1 to $2s$ to get
\begin{align*}
 sR(x,s) \ln 2s & \leq \int_1^{2s} R(x,t)dt \\
                 & =-F(x,2s)-\int_0^1 R(x,t)dt \ \ \ \ \text{ by } \eqref{eq:F=-Int_R}.
\end{align*}
Then using estimate $\eqref{eq:3}$, we have
\[
R(x,s)\leq \frac{-F(x,2s)}{s\ln2s}-\frac{\int_0^1
R(x,t)dt}{s\ln 2s}=o(1/s) \ \ \text{ as } s\rightarrow \infty,
\]
and hence estimate $\eqref{eq:4}$ follows.  \hfill
\end{proof}

\bigskip

Applying the sublemma  to the Riemann surface
$(\mathbb{R}^2,g(t))$ gives
\[
\lim_{t\rightarrow \infty} R(x,t)=0 \ \ \ \ \text{ for all } x\in
\mathbb{R}^2,
\]
which, together with the Harnack estimate, shows that $R(x,t)\equiv
0$ everywhere. This leads to a contradiction. The result follows.
\hfill
\end{proof}

\bigskip

By Lemma \ref{lemma:upper_bound}, we conclude that
\[
\lim_{r\rightarrow \infty}rk(o,r)=0
\]
provided that the circumference of the solution at infinity is
infinite. The theorem follows, since this is a contradiction of the
fact that
$$
\limsup_{r\rightarrow \infty} \ rk(o,r)>0.
$$  \hfill
\end{proof}

\section{The lower bound on $R_{\max}$} \label{sect:finite_circumference}

Since the scalar curvature of such solutions decays to
zero at spatial infinity, the scalar curvature attains its maximum at each time
slice.
By the Harnack estimate, the function $R_{\max}(t)=\max R(\cdot,t)$
is nondecreasing. Does a Type II ancient solution on a surface
satisfy $\lim_{t\rightarrow -\infty}R_{\max}(t)>0$? The main result
of this section, Theorem \ref{thm:c=2pi=>r>=4} below, gives an
affirmative answer to the noncompact case. By Theorem
\ref{thm:zero_aperture=finite_circumference}, we may assume without
loss of generality that $\mathcal{C}=2\pi$.

\begin{theorem} \label{thm:c=2pi=>r>=4}
If $(\mathbb{R}^{2},g(t))$ is a Type II ancient solution with
$\mathcal{C}=2\pi$, then we have $R_{\max}(t) \geq 4$ for all $t$.
\end{theorem}

\begin{proof} By the strong maximum principle and the fact that the
surface is non-flat, we have $R(g(t))>0$ for all $t$. Since
$\mathcal{A}=0$ and $\mathcal{C}=2\pi$, it follows from Theorem 7.11
in \cite{CC} that we have the following splitting:
$$
\left(\mathbb{R}^2,g(t),p_i\right) \rightarrow \mathbb{R} \times
\mathbb{S}^1 \text{ if } p_i \rightarrow \infty,
$$
which shows that the injectivity radius is at most $\pi$. This
implies that the supremum of the sectional curvature on the surface is
at least 1, hence that
\[
2\leq R_{\max}(g(t)) < \infty \text{ for all } t.
\]
This enables us to take a pointed limit of the sequence
$(\mathbb{R}^2,\tilde{g}_j(t),x_j)$, where the (unnormalized) metric
$\tilde{g}_j(t)$ is defined by
\[
\tilde{g}_j(t)=g(t_j+t),
\]
and points and times $(x_j,t_j)$ are chosen as in Lemma
\ref{lemma:backward_limit}. By Proposition
\ref{thm:backward_limit=cigar}, there exists a subsequence of
$(\mathbb{R}^2,\tilde{g}_j(t),x_j)$ which limits to a pointed limit
$(\mathbb{R}^2,\tilde{g}(t),O)$, which is a multiple of the cigar
soliton with
\begin{align*}
R_{\max}(\tilde{g})& =R_{\tilde{g}}(O)
=\lim_{j\rightarrow \infty}R_{\tilde{g}_j}(x_j)\\
&=\lim_{j\rightarrow \infty} R_{g(t_j)}(x_j) \in [2,\infty).
\end{align*}
(Note that the limit $\tilde{g}$ is independent of the choice of the
sequence $(x_j,t_j)$ if it satisfies $R(x_j,t_j)=R_{\max}(g(t_j))$ and
$t_j \rightarrow -\infty$.) Therefore, we obtain that
\begin{equation} \label{R>=2}
2 \leq R_{\max}(\tilde{g}(t)) < \infty.
\end{equation}

Indeed, we can improve the lower bound by the following.
\begin{lemma} \label{lemma:circumference_at_most_2pi}
The (unnormalized) backward limit $(\mathbb{R}^2,\tilde{g}(t),O)$ is
a multiple of the cigar soliton with  $\mathcal{C}_{\tilde{g}}\leq2\pi$.
 Moreover, we have $ 4 \leq R_{\max}(\tilde{g}).$
\end{lemma}

Combining this estimate with the Harnack estimate, we have
\[ 
R_{\max}(g(t)) \geq R_{\max}(\tilde{g}) \geq 4
\] 
for all $t$ as claimed. Theorem \ref{thm:c=2pi=>r>=4} follows. \hfill
\end{proof}

\bigskip

\begin{proof}[Proof of Lemma \ref{lemma:circumference_at_most_2pi}]
We first introduce some notation. Let $\bar{g}_k$ and $\bar{g}_k(t)$
denote the metric $g(-k)$ and the solution $g(t-k)$, respectively.
We realize the surface $(\mathbb{R}^{2},\bar{g}_k)$ as the graph of
a nonnegative strictly convex smooth function $ f_{k}$ as in
subsection \ref{subsect:aperture_circumeference_total_curvature}.
Let $\mathcal{I}_k$, $\mathcal{M}_k$ and $p_k$ denote the isometric
embedding, embedded hypersurface and point
$\mathcal{I}_k^{-1}((0,0,0))\in \mathbb{R}^2$, respectively. Since
the total curvature equals $2\pi$, the embedding $\mathcal{I}_k$ is
unique up to isometry. By convenient abuse of notation, denote by
$\{f \leq n \}$ and $\{ f=n \}$ the graphs of any given function $f$
over the sets $\{f \leq n \}$ and $\{ f=n \}$, and still call them
the sublevel set $\{f \leq n \}$ and level set $\{ f=n \}$,
respectively.

To show that $\mathcal{C}_{\tilde{g}} \leq2\pi$, we investigate the
pointed limit of the sequence $(\mathbb{R}^2,\bar{g}_k(t),p_k)$ as
follows.
\begin{sublemma} \label{lemma:circumference<2pi}
The sequence $(\mathbb{R}^2,\bar{g}_k(t),p_k)$ converges to a
complete pointed limit $(\mathbb{R}^2,\bar{g}(t),p)$ with
$\mathcal{C}_{\bar{g}} \leq 2\pi$. Consequently, we have $
\mathcal{A}_{\bar{g}}=0 \text{ and } \tau_{\bar{g}}=2\pi. $
Furthermore, the embedding $\mathcal{I}:(\mathbb{R}^2,\bar{g}(0))
\hookrightarrow \mathbb{R}^3$ is unique up to congruence.
\end{sublemma}

\begin{proof} Since $(\mathbb{R}^2,g(t))$ is a Type II ancient
solution, the complete pointed surfaces
$(\mathbb{R}^{2},\bar{g}_k(t),p_k)$ have uniformly bounded curvature
on any given finite interval containing  $t=0$\ and satisfy the
injectivity radius estimate at $t=0$. By Hamilton's compactness
theorem, the sequence
$(\mathbb{R}^{2},\bar{g}_k(t),p_k)$ subconverges to an eternal
solution $(\mathbb{R}^{2},\bar{g}(t),p)$ with uniformly bounded
curvature. Denote by $\bar{g}$ the metric $\bar{g}(0)$. As noted,
the surface $(\mathbb{R}^{2},\bar{g})$ is realizable  as the graph
of a nonnegative strictly convex smooth function $f$. Let
$\mathcal{I}$ and $\mathcal{M}$ denote the isometric embedding and
embedded hypersurface, respectively. It is clear that we have
$\mathcal{I}(p)=(0,0,0)\in \mathcal{M} \subset \mathbb{R}^3$.

By the strong maximum principle we have either $R_{\bar{g}(t)}\equiv
0$ or $R_{\bar{g}(t)}>0$ everywhere in space-time. We now claim that
the curvature is positive. The sequence
$(\mathcal{M}_{k},\bar{g}_k,O)$ subconverges to the corresponding
embedded pointed surface $(\mathcal{M},\bar{g},O)$. Note that if
$R_{\bar{g}}\equiv 0$, then the embedded surface must be a cylinder
or a plane, which is impossible because the surface $\mathcal{M}_k$
is always inside the set $\{ (x,y,z)\in \mathbb{R}^3\ | \ x^2+y^2
\leq \pi^2,\ z \geq 0 \}$. Therefore, $R_{\bar{g}}>0$. In other
words, the solution $(\mathbb{R}^{2},\bar{g})$ is non-flat as
claimed.

Since the graph of $f_{k}$ subconverges to the graph of $f$, the
level curve $f_{k}^{-1}\{c\}$ uniformly subconverges to the level
curve $f^{-1}\{c\}$. As seen in subsection
\ref{subsect:aperture_circumeference_total_curvature}, the length of
level curves of $f$ converges to $\mathcal{C}_{\bar{g}},$ thus for
any $\eta >0$ there exists a positive constant $c_0$ such that
\[
\ell(\{ f =c \}) > \mathcal{C}_{\bar{g}}-\eta
\]
for all $c \geq c_0$. On the other hand, by Theorem 1 in \cite{AZ}
(p225), we have
\[
 \liminf_{k\rightarrow \infty} \ell(\{ f_k =c \}) \geq \ell(\{ f =c
 \}),
\]
thus there also exits a positive integer $k_0$ such that
\[
\ell(\{ f_{k_0} =c \}) > \mathcal{C}_{\bar{g}}-2\eta
\]
for all $c>c_0$. This implies that
\begin{align*}
2\pi  =\mathcal{C}_{\bar{g}_{k_0}}  & > \ell(\{ f_{k_0} =c \}) \\
& > \mathcal{C}_{\bar{g}}-2\eta
\end{align*}
for all $\eta>0$. Therefore, we see that $\mathcal{C}_{\bar{g}} \leq
2\pi$ as claimed.  The rest of the proof is immediate from
Lemma \ref{prop:embbedding_rigid}. \hfill
\end{proof}

\medskip

Intuitively, the point $p_k$ should be close to the point, denoted
by $q_k$, where the curvature $R_{\bar{g}_k}$ attains its maximum
at. This motivates us to take a pointed limit of the sequence
$(\mathbb{R}^{2},\bar{g}_k,q_k)$ and have the following.

\begin{sublemma} \label{lemma:last_one}
The sequence $(\mathbb{R}^{2},\bar{g}_k,q_k)$ converges to a pointed
limit $(\mathbb{R}^{2},\bar{g},\bar{q})$.
\end{sublemma}

\begin{proof} We first claim that the set $\{p_{k,\text{ }}p_{j,}$
$q_{j}\}_{j= k}^\infty$ is uniformly bounded in the pointed surface
$(\mathbb{R}^{2},\bar{g}_k,p_{k})$, $k\geq N$, for some positive
integer $N$.

It follows from Sublemma \ref{lemma:circumference<2pi} that for any
given $\varepsilon>0$, there exist a positive number $d>0$ and a
positive integer $N$ such that
\begin{subequations}
\begin{equation} \label{eq:almost_2pi_1}
\int_{\{f<d\}}K_{\bar{g}}dv_{\bar{g}}>2\pi-\varepsilon
\end{equation}
and
\begin{equation} \label{eq:almost_2pi_2}
\int_{\{f_{k}<d\}}K_{\bar{g}_k}dv_{\bar{g}_k}>2\pi-2\varepsilon
\end{equation}
\end{subequations}
for $k\geq N$.

Since $R_{\bar{g}_k}(p_{k})\rightarrow R_{\bar{g}}(O)>0$ as $k
\rightarrow \infty$, there exists a positive integer $k_1 \geq N$
such that
\begin{equation}\label{eq:almost_1/2_R}
R_{\bar{g}_k}(p_{k}) >\frac{1}{2} R_{\bar{g}}(O)
\end{equation}
for $k\geq k_1$. Combining \eqref{eq:almost_1/2_R}, \eqref{R>=2} and
the Harnack estimate \cite{H_harnack}, we get the following
estimates:
\begin{subequations}
\begin{equation}\label{eq:R_g_bar>=2}
R_{\bar{g}_k}(q_j) \geq R_{\bar{g}_j}(q_j)=R_{\max}(\bar{g}_j)\geq 2
\end{equation}
and
\begin{equation}\label{R_g_k_bar>=1/2R_g_bar}
R_{\bar{g}_k}(p_{j}) \geq R_{\bar{g}_j}(p_{j}) >\frac{1}{2}
R_{\bar{g}}(O)
\end{equation}
\end{subequations}
for all $j \geq k \geq k_1$. As a consequence of the injectivity
radius control on $\mathcal{M}_{k}$,  gradient estimate on scalar
curvature $R_{\bar{g}_k}$ and \eqref{eq:almost_2pi_1},
\eqref{eq:almost_2pi_2}, \eqref{eq:R_g_bar>=2},
\eqref{R_g_k_bar>=1/2R_g_bar},  we obtain uniform estimates on
$$
\mathrm{dist}(p_{j},\{ f_{k}<d \}) \text{ and } \mathrm{dist}(q_{j},
\{ f_{k}<d \})
$$
for all $j\geq k \geq k_1$, which implies that the set
$\{p_{k,\text{ }}p_{j,}$ $q_{j}\}_{j= k}^\infty$ is uniformly
bounded in the pointed surface $(\mathbb{R}^{2},\bar{g}_k,p_{k})$,
$k\geq k_1$. Therefore, the sequence
$(\mathbb{R}^{2},\bar{g}_k,q_k)$ converges to a pointed limit
$(\mathbb{R}^{2},\bar{g},\bar{q})$. \hfill
\end{proof}

\medskip

Note that, as pointed out, the sequence
$(\mathbb{R}^{2},\bar{g}_k,q_k)$ converges to the pointed limit
$(\mathbb{R}^{2},\tilde{g},p)$. By the rotational symmetry of the
cigar soliton, we have $\bar{q}=p$ and hence $\tilde{g}$ coincides
with $\bar{g}$. From Sublemma \ref{lemma:circumference<2pi}, the
circumference of $(\mathbb{R}^2,\tilde{g}(t))$ at infinity is at
most $2\pi$ as claimed. Since $(\mathbb{R}^2,\tilde{g})$ is a
multiple of the cigar soliton with $\mathcal{C}_{\tilde{g}}\leq
2\pi$, Lemma \ref{lemma:circumference_at_most_2pi} follows
immediately. \hfill
\end{proof}

\medskip

As a result of Theorem \ref{thm:c=2pi=>r>=4}, the scalar curvature
assumes its maximum in space-time provided that $\mathcal{C}=2\pi$ and
$R(x,t) \leq 4$. This means that the solution is indeed a
Type II singularity model. By Hamilton's theorem \cite{H_eternal},
such a solution must be the cigar soliton. Therefore, we obtain the
following.
\begin{corollary} \label{cor:c=2pi&r=<4=>soliton}
If $(\mathbb{R}^{2},g(t))$ is a  Type
II ancient solution with $\mathcal{C}=2\pi$ and $R(x,t) \leq 4$
 for all $(x,t)\in (-\infty,\infty) \times
\mathbb{R}^2$,  then it is the cigar soliton.
\end{corollary}

\medskip

We end this section by the following.

\medskip

\noindent\textbf{Remark}. S. Angenent and L. Wu  \cite{W} observe
that the logarithmic fast diffusion equation
\begin{equation}
\dfrac{\partial }{\partial t}u= \triangle \log u
\label{eq:log-fast-diffusion}
\end{equation}
on the plane $\mathbb{R}^2$, where $\triangle$ denotes the Euclidean
Laplace operator on $\mathbb{R}^2$, represents the evolution of the
conformally flat metric with
$$ g=u(dx^2+dy^2)$$
under the Ricci flow. The equivalence follows from the facts that
the conformal metric $g$ has scalar curvature $R =-(\triangle\log
u)/u$ and in two dimensions $R_{ij} = \dfrac{1}{2} Rg_{ij}.$
Daskalopoulos and Sesum \cite{DS} study the classification of
eternal solutions of equation \eqref{eq:log-fast-diffusion}. They
show that any positive smooth eternal solution $u(x,y,t)$ is a
gradient soliton of the form
$$u(x,y,t)=\dfrac{2}{\beta(|x-x_0|^2+|y-y_0|^2+\delta e^{2\beta t})} $$
for some $(x_0,y_0) \in \mathbb{R}^2$ and some positive constants
$\beta, \delta$, provided that the solution $u$ defines a complete
metric of bounded curvature and bounded width. Note that Theorem
\ref{thm:zero_aperture=finite_circumference} removes  the hypothesis
of the  width  being finite.

\medskip

\noindent \textbf{Acknowledgements}. The author would like to thank
Professor Ben Chow for bringing the idea of taking limits backwards
in time, Ni's work \cite{Ni} and Ni and Tam's notes \cite{NT2,NT3}  to
his attention. The author also would like to thank Professor L. Ni
and Professor L.F. Tam for their indispensable notes, and thank
Professor Peng Lu for his helpful suggestions and comments. Finally,
the author thanks the referee for the helpful suggestions concerning the
presentation of this paper.


\begin{thebibliography}{999}

\bibitem[AZ]{AZ} A.D. Alexandrov and V.A. Zalgaller, \textit{Intrinsic
geometry of surfaces}, Volume \textbf{15}, Translations of
Mathematical Monographs.

\bibitem[BZ]{BZ} Y.D. Burago and V.A. Zalgaller, \textit{Geometric
inequality}, Springer Series in Soviet Mathematics,
Springer-Verlag, Berlin, 1988.

\bibitem[CC]{CC} J. Cheeger and T. Colding, \textit{Lower bounds on Ricci
curvature and the almost rigidity of warped products}, Ann. of
Math. \textbf{144} (1996), 189-237.

\bibitem[CLN]{CLN} B. Chow, P. Lu and L. Ni, \textit{Hamilton's Ricci flow},
to appear.

\bibitem[DH]{DH} P. Daskalopoulo and R. Hamilton,
\textit{Geometric estimates for the logarithmic fast diffusion
Equation}, Comm. Anal.  Geom. \textbf{12} (2004), 143-164.

\bibitem[DS]{DS} P. Daskalopoulos and N. Sesum, \textit{Eternal solutions
to the Ricci flow on} $\mathbb{R}^2$, arXiv:math.AP/0603525.

\bibitem[GS]{GS}R. Greene and K. Shiohama, \textit{Convex functions on
complete noncompact manifolds: Topological structure}, Invent. Math.
\textbf{63} (1981), 129-157.

\bibitem[H1]{H_surface} R. Hamilton, \textit{The Ricci flow on surfaces},
Contemporary Mathematics \textbf{71} (1988), 237-261.

\bibitem[H2]{H_harnack} R. Hamilton, \textit{The Harnack estimate for the
Ricci flow}, J. Diff. Geom. \textbf{37} (1993), 225-243.

\bibitem[H3]{H_eternal} R. Hamilton, \textit{Eternal solutions to the Ricci
flow}, J. Diff. Geom. \textbf{38} (1993), 1-11.

\bibitem[H4]{H_sing} R. Hamilton, \textit{The formation of singularities in
the Ricci flow}, Surveys in Differential Geometry \textbf{2} (1995),
7-136, International Press.

\bibitem[H5]{Hartman} P. Hartman, \textit{Geodesic parallel coordinates in
the large}, Amer. J. Math. \textbf{86} (1964), 705-727.

\bibitem[Ni]{Ni} L. Ni, \textit{Ancient solutions to K\"{a}hler-Ricci flow},
Math. Res. Lett. \textbf{12} (2005), 633-653.

\bibitem[NT1]{NT1} L. Ni and L.F. Tam, \textit{K\"{a}hler-Ricci flow and the
Poincar\'{e}-Lelong equation}, Comm. Anal. Geom. \textbf{12} (2004),
111-141.

\bibitem[NT2]{NT2} L. Ni and L.F. Tam, \textit{Notes on soliton} 1,
unpulished.

\bibitem[NT3]{NT3} L. Ni and L.F. Tam, \textit{Notes on soliton} 2,
unpulished.

\bibitem[O]{O} S.P. Olovjani\v{s}nikov, \textit{On the bending of infinite convex surfaces},
Mat. Sb. \textbf{18} (\textbf{60}) (1946), 429-440.

\bibitem[P]{P} A.V. Pogorelov, \textit{Extrinsic geometry of convex surface},
Volume \textbf{35}, Translations of Mathematical Monographs.

\bibitem[R]{R} P. Rosenau, \textit{On fast and super-fast diffusion}, Phys.
Rev. Lett. \textbf{74} (1995), 1056-1059.

\bibitem[Sh]{Sh} W.-X. Shi, \textit{Ricci flow and the uniformization on
complete noncompact K\"{a}hler manifolds}, J. Diff. Geom.
\textbf{45} (1997), 94-220.

\bibitem[St]{St}J.J. Stoker, \textit{\"{U}ber die Gestalt der positiv
gekr\"{u}mmten offenen Fl\"{a}chen im dreidimensionalen Raum}, Comp.
Math. \textbf{3} (1936), 55-88.

\bibitem[W]{W} L.F. Wu, \textit{The Ricci flow on complete $\mathbb{R}^{2}$},
Comm. Anal. Geom. \textbf{1},  (1993), 439-472.
\end{thebibliography}
\end{document}